\documentclass[a4paper,12pt]{amsart}
\usepackage{amssymb,amscd}

\calclayout

\theoremstyle{plain}
\newtheorem{theo}{Theorem}[section]
\newtheorem{prop}[theo]{Proposition}
\newtheorem{lemm}[theo]{Lemma}
\newtheorem{coro}[theo]{Corollary}

\theoremstyle{definition}

\theoremstyle{remark}
\newtheorem*{rema}{Remark}
\newtheorem*{clai}{Claim}

\newcommand{\field}[1]{\mathbb{#1}}
\newcommand{\C}{\field{C}}
\newcommand{\F}{\field{F}}
\newcommand{\Q}{\field{H}}
\newcommand{\R}{\field{R}}
\newcommand{\Z}{\field{Z}}
\newcommand{\CS}{\mathcal{J}}

\newcommand{\side}[1]{{}^{#1}\hskip-.075em}
\newcommand{\siden}[1]{{}^{#1}\hskip-.375em}

\DeclareMathOperator{\GL}{GL} \DeclareMathOperator{\Hom}{Hom}
\DeclareMathOperator{\Irr}{Irr} \DeclareMathOperator{\Vect}{Vect}
\DeclareMathOperator{\ind}{ind} \DeclareMathOperator{\res}{res}
\DeclareMathOperator{\Rep}{Rep} \DeclareMathOperator{\diag}{diag}

\begin{document}

\title{Classification of Equivariant Real Vector Bundles over a Circle}

\date{May 8, 2000}

\author{Jin-Hwan Cho}
\address{School of Mathematics, Korea Institute for Advanced Study}
\email{chofchof@kias.re.kr}
\thanks{Jin-Hwan Cho was supported by postdoctoral fellowships program
from Korea Science \& Engineering Foundation (KOSEF)}

\author{Sung Sook Kim}
\address{Department of Applied Mathematics, Paichai University}
\email{sskim@www.paichai.ac.kr}

\author{Mikiya Masuda}
\address{Department of Mathematics, Osaka City University}
\email{masuda@sci.osaka-cu.ac.jp}

\author{Dong Youp Suh}
\address{Department of Mathematics, Korea Advanced Institute of
Science and Technology}
\email{dysuh@math.kaist.ac.kr}
\thanks{Dong Youp Suh was partially supported from the inter disciplinary
Research program of the KOSEF grant no.~(1999-2-101-002-4)}

\subjclass{Primary 57S25; Secondary 19L47}

\keywords{group action, equivariant real vector bundle, circle,
fiber module, extension of representation}

\begin{abstract}
This is a continuation of the authors' previous
work~\cite{CKMS99a} on classification of equivariant complex
vector bundles over a circle. In this paper we classify
equivariant real vector bundles over a circle with a compact Lie
group action, by characterizing the fiber representations of
them, and by using the result of the complex case. We also treat
the triviality of them. The basic phenomenon is similar to the
complex case but more complicated here.
\end{abstract}

\maketitle

\section{Introduction}

In~\cite{CKMS99a} we classified equivariant \emph{complex} vector
bundles over a circle, and in this paper we classify equivariant
\emph{real} ones. The argument developed in this paper is similar
to that in~\cite{CKMS99a} but is rather more complicated. The
complexity arises from two aspects: one is topology and the other
is representation theory. For instance, any (nonequivariant)
complex line bundle over a circle is trivial while there are two
non-isomorphic real line bundles, that is, the Hopf line bundle
and the trivial one. This is an evidence of the topological
complexity in the real case. As for representation theory, it is
recognized in general that real representation theory is more
complicated than complex representation theory.

Let us introduce some notation to state our results. Let $G$ be a
compact Lie group and let $\rho\colon G\to O(2)$ be an orthogonal
representation. The unit circle of the corresponding $G$-module is
denoted by $S(\rho)$. It is well known that any circle with
$G$-action is equivalent to $S(\rho)$ for some $\rho$. We set
$H=\rho^{-1}(1)$, so that $H$ acts trivially on $S(\rho)$ and the
\emph{fiber $H$-module} of a real $G$-vector bundle over $S(\rho)$
is determined uniquely up to isomorphism.

Let $\Irr(H)$ be the set of characters of irreducible real
$H$-modules. It has a $G$-action defined as follows: For $\chi\in
\Irr(H)$ and $g\in G$, $\side{g}\chi\in\Irr(H)$ is defined by
$\side{g}\chi(h)=\chi(g^{-1}hg)$ for $h\in H$. Since a character
is a class function, the isotropy subgroup $G_\chi$ of $G$ at
$\chi\in\Irr(H)$ contains $H$. We choose and fix a representative
from each $G$-orbit in $\Irr(H)$ and denote the set of those
representatives by $\Irr(H)/G$. Denote by $\Vect_G(S(\rho))$ the
set of isomorphism classes of real $G$-vector bundles over
$S(\rho)$ and by $\Vect_{G_{\chi}}(S(\rho),\chi)$ the subset of
$\Vect_{G_\chi}(S(\rho))$ with a multiple of $\chi$ as the
character of fiber $H$-modules. They are semi-groups under Whitney
sum. The decomposition of a $G$-vector bundle into the
$\chi$-isotypical components induces an isomorphism
\[
\Vect_G(S(\rho)) \cong
\prod_{\chi\in\Irr(H)/G}\Vect_{G_{\chi}}(S(\rho),\chi),
\]
see~\cite[Section~2]{CKMS99a}. This reduces the study of
$\Vect_G(S(\rho))$ to that of $\Vect_{G_\chi}(S(\rho),\chi)$, and
since $\chi$ is $G_\chi$-invariant and $G_\chi$ is again a
compact Lie group, we are led to study $\Vect_G(S(\rho),\chi)$
where $\chi$ is $G$-invariant, namely $\side{g}\chi=\chi$ for all
$g\in G$.

Let $E\in\Vect_G(S(\rho),\chi)$. The fiber $H$-module of $E$ has a
multiple of $\chi$ as the character by definition. In fact, the
fiber $E_z$ of $E$ over a point $z\in S(\rho)$ is a real module of
the isotropy subgroup $G_z$ at $z$, and unless $\rho(G)\subset
SO(2)$, $G_z$ properly contains $H$ for some $z$. It turns out
that these fiber $G_z$-modules almost distinguish elements in
$\Vect_G(S(\rho),\chi)$. To be more specific, we shall introduce
some more notation. Unless $\rho(G)\subset SO(2)$, $\rho(G)$ is
$O(2)$ or a dihedral group $D_n$ of order $2n$ for some
positive integer $n$. We identify $S(\rho)$ with the unit circle
of the complex line $\C$, so that the dihedral group $D_n$ is
generated by the rotation through an angle $2\pi/n$ and the
reflection about the $x$-axis. Then the isotropy subgroups $G_z$
at $z=1$ (and $z=e^{\pi i/n}$ when $\rho(G)=D_n$) contain $H$ as
an index two subgroup unless $\rho(G)\subset SO(2)$. For a group
$K$ containing $H$ we denote by $\Rep(K,\chi)$ the set of
isomorphism classes of real $K$-modules whose characters
restricted to $H$ are multiples of $\chi$. The set $\Rep(K,\chi)$
is a semi-group under direct sum. Restriction of elements in
$\Vect_G(S(\rho),\chi)$ to fibers at $1$ (and $\mu=e^{\pi i/n}$
when $\rho(G)=D_n$) yields a semi-group homomorphism
\[
\Gamma\colon \Vect_G(S(\rho),\chi) \to
\begin{cases}
\Rep(H,\chi), &\quad\text{if $\rho(G)\subset SO(2)$}, \\
\Rep(G_1,\chi), &\quad\text{if $\rho(G)=O(2)$}, \\
\Rep(G_1,G_\mu,\chi), &\quad\text{if $\rho(G)=D_n$},
\end{cases}
\]
where $\Rep(G_1,G_\mu,\chi)$ denotes the subsemi-group of
$\Rep(G_1,\chi)\times \Rep(G_\mu,\chi)$ consisting of pairs of the
same dimension. The map $\Gamma$ can also be defined in the
complex case and is proved to be an isomorphism
in~\cite[Proposition~6.2]{CKMS99a}. However, $\Gamma$ is not
always an isomorphism in the real case, for instance the two
non-isomorphic line bundles over a circle (when $G$ is trivial)
mentioned before have obviously the same $\Gamma$ image.
Nevertheless it turns out that $\Gamma$ is an isomorphism in most
cases. Remember that $\chi$ is called of \emph{real},
\emph{complex}, or \emph{quaternion} type if the $H$-endomorphism
algebra of the irreducible real $H$-module with $\chi$ as the
character is isomorphic to $\R$, $\C$, or $\Q$ respectively, and
that $\chi$ is called $K$-\emph{extendible} if it extends to a
character of a group $K$ containing $H$. There are two cases in
which the classification works somewhat exceptionally.

\smallskip

\begin{flushleft}
\emph{Case A:} $\rho(G)\subsetneq SO(2)$ and $\chi$ is of real type. \\
\emph{Case B:} $\rho(G)=D_n$, $\chi$ is of real type and
neither $G_1$- nor $G_\mu$-extendible.
\end{flushleft}

\smallskip

\begin{theo} \label{theo:main_theorem}
Except for Cases~A and~B, the semi-group homomorphism $\Gamma$ is
an isomorphism. In Case~A or~B, $\Gamma$ is a two to one map;
more precisely, there is a free involution on
$\Vect_G(S(\rho),\chi)$ given by tensoring with a nontrivial
$G$-line bundle (with trivial fiber $H$-module) and $\Gamma$
induces an isomorphism on the orbit space.
\end{theo}

Theorem~\ref{theo:main_theorem} reduces the study of
$\Vect_G(S(\rho),\chi)$ to that of real representation theory,
especially to the study of $\Rep(K,\chi)$ where $K$ is a compact
Lie group containing $H$ as an index two subgroup and $\chi$ is
$K$-invariant. The complexity of real representation theory
emerges here. Namely, the number of $K$-extensions of $\chi$ can
be zero, one, or two, while it is always two in the complex case.
Combining this observation with Theorem~\ref{theo:main_theorem},
one sees that the semi-group structures on $\Vect_G(S(\rho),\chi)$
are of five types depending on $\rho(G)$ and $\chi$ (see
Theorem~\ref{theo:structure}) while they are of three types in
the complex case (see~\cite[Theorem~B]{CKMS99a}).

The paper is organized as follows. In Section~2 we shall determine
the semi-group structure of the target space of the semi-group
homomorphism $\Gamma$. For this we need to study some
representation theory, especially on extensions of
representations. The semi-group structure of the target space of
$\Gamma$ is given in Lemma~\ref{lemm:structure_Rep(K)} and
Lemma~\ref{lemm:structure_Rep(K1,K2)}.

In section~3 we prove Theorem~\ref{theo:main_theorem} except for
Cases~A and~B. Indeed, Proposition~\ref{prop:construction} shows
that $\Gamma$ is always surjective, and
Proposition~\ref{prop:real_isomorphism_theorem} shows that
$\Gamma$ is injective except for Cases~A and~B.

The proof for Cases~A and~B are treated in Section~4. Case~A can
be easily proved by reducing the case to nonequivariant case. For
Case~B we consider possible complex $G$ vector bundle structures
on an element of $\Vect_G(S(\rho),\chi)$. Then we use the
classification results of complex $G$ vector bundles over a circle
in \cite{CKMS99a}, and some counting argument to finish the proof
for Case~B.

The semi-group structure of $\Vect_G(S(\rho),\chi)$ is given in
Theorem~\ref{theo:structure} of Section~5 by figuring out the
generators with their relations.
Theorem~\ref{theo:triviality_for_L} in Section~6 shows which of
the generators of  $\Vect_G(\S(\rho),\chi)$ are trivial.

\section{Ingredients from representation theory}

We shall determine the semi-group structure of the target space
of $\Gamma$. For that, the following lemma from representation
theory plays a key role. Throughout this section, $H$ is an index
two normal subgroup of a group $K$ and $U$ is a real irreducible
$H$-module with $K$-invariant character. By the type of $U$ we
mean the type of the character of $U$. As is well known, any
$K$-extension of $U$ appears in $\ind_H^KU$ as a direct summand
at least once, which follows from the Frobenius reciprocity, and
$\res_H\ind_H^KU\cong 2U$ because $H$ is of index two and the
character of $U$ is $K$-invariant.

\begin{lemm} \label{lemm:real_Mackey}
\textnormal{(1)} Suppose $U$ is $K$-extendible.
\begin{enumerate} \renewcommand{\theenumi}{\alph{enumi}}
\item If $U$ is of real type, there are two mutually
non-isomorphic $K$-extensions of real type.
\item If $U$ is of complex type, there are either two mutually
non-isomorphic $K$-extensions of complex type or unique
$K$-extension of real type, not both.
\item If $U$ is of quaternionic type, there are either two
mutually non-isomorphic $K$-extensions of quaternionic type or
unique $K$-extension of complex type, not both.
\end{enumerate}
In any case, if $U$ has two $K$-extensions, then one is isomorphic to the
tensor product of the other with the nontrivial real $K$-module of dimension
one with trivial $H$-action.

\textnormal{(2)} Suppose $U$ is not $K$-extendible. Then $U$ is
not of quaternionic type, and $\ind_H^K U$ is irreducible of
complex or quaternionic type according as $U$ is of real or
complex type, respectively. Moreover, any $K$-extension of $2U$
is isomorphic to $\ind_H^K U$.
\end{lemm}

\begin{proof}
(1) Let $W$ be a $K$-extension of $U$. Since $U$ is irreducible,
so is $W$. Set $d(W)=\dim_\R\Hom_K(W,W)$ and
$d(U)=\dim_\R\Hom_H(U,U)$. Note that $d(U)=1,2$, or $4$ according
as $U$ is of real, complex, or quaternionic type, respectively.
Since $\res_HW\cong U$, it follows from the Frobenius reciprocity
that
\begin{equation}
m(W)\cdot d(W)=\dim_\R\Hom_K(W,\ind_H^K
U)=\dim_\R\Hom_H(\res_HW,U)=d(U), \tag{*}
\end{equation}
where $m(W)$ denotes the multiplicity of $W$ in $\ind_H^K U$.
On the other hand, since
\[
\dim_\R\ind_H^K U=2\dim_\R U=2\dim_\R W,
\]
the multiplicity $m(W)$ must be either $1$ or $2$.

If $m(W)=1$ (i.e., $d(W)=d(U)$), then $\ind_H^K U\cong W\oplus W'$
for some $K$-module $W'$ which is not isomorphic to $W$. Since
$\res_H\ind_H^K U\cong 2U$, $\res_HW\cong U$ is isomorphic to
$\res_HW'$, so that $W'$ is also a $G$-extension of $U$. Since
$W'$ appears in $\ind_H^KU$ once, $m(W')=1$.  Therefore the
identity (*) applied to $W'$ implies $d(W')=d(U)$. This together
with the equality $d(W)=d(U)$ implies that the two mutually
non-isomorphic $K$-extensions $W$ and $W'$ of $U$ are of the same
type as $U$.

If $m(W)=2$, then $\ind_H^K U\cong 2W$. Since any $K$-extension of
$U$ is contained in $\ind_H^K U$ as a direct summand, $W$ is the
unique $K$-extension of $U$. The type of $W$ can be read from the
equality $d(W)=d(U)/2$. This equality in particular implies that
$d(U)$ must be even, in other words, $U$ is not of real type when
$m(W)=2$.

These prove the statement~(1) except the last statement.  The
last statement can be seen as follows. If $U$ has two
$K$-extensions $W$ and $W'$, then $\Hom_H(W,W')$ is a nontrivial
real $K$-module of dimension one with trivial $H$-action, and
$W\otimes\Hom_H(W,W')$ is isomorphic to $W'$, proving the last
statement.

(2) Suppose $U$ has no $K$-extension. It is known that $U$ is
$K$-extendible if $U$ is of quaternionic type, see the remark
following Proposition~4.5 in~\cite{CKS99}. If $\ind_H^K U$ is
reducible, then each direct summand of $\ind_H^K U$ is a
$K$-extension of $U$ which contradicts the assumption that $U$
has no $K$-extension. Therefore, $\ind_H^K U$ is irreducible.
Noting that $\res_H\ind_H^K U\cong 2U$, it follows from the
Frobenius reciprocity that
\[
d(\ind_H^K U)=\dim_\R\Hom_K(\ind_H^K U,\ind_H^K
U)=\dim_\R\Hom_H(2U,U)=2d(U),
\]
which implies the statement on the type of $\ind_H^K U$. The last
statement in~(2) follows again from the Frobenius reciprocity.
\end{proof}

The following lemma follows from Lemma~\ref{lemm:real_Mackey}.

\begin{lemm} \label{lemm:structure_Rep(K)}
Let $\chi$ be the $K$-invariant character of $U$, and let $e$ be
the number of $K$-extensions of $U$. Then $e=0$, $1$, or $2$, and
\begin{enumerate}
\item if $e=0$, then $\Rep(K,\chi)$ is generated by $\ind_H^K U$,
\item if $e=1$, then $\Rep(K,\chi)$ is generated by a $K$-module of the
same dimension as $U$,
\item if $e=2$, then $\Rep(K,\chi)$ is generated by two $K$-modules
of the same dimension as $U$ such that one is isomorphic to the
tensor product of the other with the nontrivial real $K$-module
of dimension one with trivial $H$-action. \qed
\end{enumerate}
\end{lemm}

Let $K_1$ and $K_2$ be two groups containing $H$ as an index two
subgroup, and let $\chi$ be a real irreducible character of $H$
which is $K_1$- and $K_2$-invariant. We next consider the subset
$\Rep(K_1,K_2,\chi)$ of $\Rep(K_1,\chi)\times\Rep(K_2,\chi)$
consisting of pairs of the same dimension. Denote by $e_1$ and
$e_2$ the number of $K_1$- and $K_2$-extensions,
respectively.

\begin{lemm} \label{lemm:structure_Rep(K1,K2)}
The semi-group $\Rep(K_1,K_2,\chi)$ is generated by
\[
\begin{cases}
\text{one element $R_\chi$}, & \text{if
$(e_1,e_2)=(0,0),(1,0),(0,1)$, or $(1,1)$}, \\
\text{two elements $R_\chi^\pm$}, & \text{if $(e_1,e_2)=(2,1)$ or
$(1,2)$}, \\
\text{three elements $\widetilde R_\chi^0,\widetilde
R_\chi^\pm$}, & \text{if $(e_1,e_2)=(2,0)$ or $(0,2)$}, \\
\text{four elements $R_\chi^{\pm\pm}$}, & \text{if $(e_1,e_2)=(2,2)$},
\end{cases}
\]
with relations $2\widetilde R_\chi^0=\widetilde R_\chi^++\widetilde
R_\chi^-$ and $R_\chi^{++}+R_\chi^{--}=R_\chi^{+-}+R_\chi^{-+}$.
\end{lemm}

\begin{proof}
For $i=1$ and $2$, denote by $\widetilde R_i$, $R_i$, and
$R_i^\pm$ the set of generators of $\Rep(K_i,\chi)$ according to
$e_i=0$, $1$, and $2$, respectively. Note that the dimension of
$\widetilde R_i$ is twice that of $R_i$ and $R_i^\pm$. Then it is easy to
find the generators and relations of $\Rep(K_1,K_2,\chi)$
according to $(e_1,e_2)$ as in Table~\ref{table:generators}.
\end{proof}

\begin{table}[hbt]
\begin{tabular}{|c|l|c|} \hline
$(e_1,e_2)$ & \multicolumn{1}{|c|}{\rule[-1ex]{0pt}{3ex}
Generators} & \multicolumn{1}{|c|}{Relations} \\\hline\hline
$(0,0)$ & \rule[-1.25ex]{0pt}{4ex} $R_\chi=(\widetilde
R_1,\widetilde R_2)$ & \\\cline{1-2} $(1,0)$ &
\rule[-1.25ex]{0pt}{4ex} $R_\chi=(R_1\oplus R_1,\widetilde R_2)$ &
\\\cline{1-2} $(0,1)$ & \rule[-1.25ex]{0pt}{4ex}
$R_\chi=(\widetilde R_1,R_2\oplus R_2)$ &
\raisebox{2ex}[0pt]{none} \\\cline{1-2} $(1,1)$ &
\rule[-1.25ex]{0pt}{4ex} $R_\chi=(R_1,R_2)$ & \\\hline\hline
$(2,1)$ & \rule[-1.25ex]{0pt}{4ex} $R_\chi^\pm=(R_1^\pm,R_2)$ &
\\\cline{1-2} $(1,2)$ & \rule[-1.25ex]{0pt}{4ex}
$R_\chi^\pm=(R_1,R_2^\pm)$ & \raisebox{2ex}[0pt]{none}
\\\hline\hline &
\rule[-1.25ex]{0pt}{4ex} $\widetilde R_\chi^0=(R_1^+\oplus
R_1^-,\widetilde R_2)$, & \\ \raisebox{2ex}[0pt]{$(2,0)$} &
\rule[-1.25ex]{0pt}{4ex} $\widetilde R_\chi^\pm=(R_1^\pm\oplus
R_1^\pm,\widetilde R_2)$ & \\\cline{1-2} &
\rule[-1.25ex]{0pt}{4ex} $\widetilde R_\chi^0=(\widetilde
R_1,R_2^+\oplus R_2^-)$, & \raisebox{2ex}[0pt]{$2\widetilde
R_\chi^0=\widetilde R_\chi^++\widetilde R_\chi^-$} \\
\raisebox{2ex}[0pt]{$(0,2)$} & \rule[-1.25ex]{0pt}{4ex}
$\widetilde R_\chi^\pm=(\widetilde R_1,R_2^\pm\oplus R_2^\pm)$ &
\\\hline\hline $(2,2)$ & \rule[-1.25ex]{0pt}{4ex}
$R_\chi^{\pm\pm}=(R_1^\pm,R_2^\pm)$ &
$R_\chi^{++}+R_\chi^{--}=R_\chi^{+-}+R_\chi^{-+}$
\\\hline
\end{tabular} \vspace{1.25ex}
\caption{Generators and relations of $\Rep(K_1,K_2,\chi)$
according to $(e_1,e_2)$} \label{table:generators}
\end{table}

\begin{rema}
For $i=1$ and $2$, denote by $\R_i^+$ and $\R_i^-$, respectively,
the trivial and the nontrivial real $K_i$-module of dimension one
with trivial $H$-action. Then the set of pairs
$(\R_1^\pm,\R_2^\pm)$ forms a group isomorphic to $\Z/2\times\Z/2$
under tensor product on each factor, and it acts by the same
operation on the generators in
Lemma~\ref{lemm:structure_Rep(K1,K2)}. The action is transitive
except for the third case where $\Rep(K_1,K_2,\chi)$ is generated
by three elements. In that case, $\widetilde R_\chi^\pm$
constitute one orbit and $\widetilde R_\chi^0$ is fixed by the
action of the pairs $(\R_1^\pm,\R_2^\pm)$.
\end{rema}

We recall some facts on the extension of an $H$-module, which will
be used in Section~6.

\begin{lemm} \label{lemm:extension}
Let $H$ be a normal subgroup of $G$ and let $U$ be a real
irreducible $H$-module with $G$-invariant character.
\begin{enumerate}
\item Suppose $G/H$ is finite cyclic of odd order. Then
$U$ has a $G$-extension, and the $G$-extension is unique if $U$
is of real type.
\item Suppose $G/H$ is a dihedral group of order $2n$ for odd $n$.
Then
\begin{enumerate}
\item $U$ has a $G$-extension if and only if it has a
$K$-extension for some subgroup $K$ of $G$ which contain $H$ as
an index two subgroup,
\item $2U$ always has a $G$-extension.
\end{enumerate}
\end{enumerate}
\end{lemm}

\begin{proof}
See~\cite[Proposition~4.5 and~4.6]{CKS99} for the former
statement in~(1) and~(2-a).

To see the uniqueness in~(1), we note that if $U_1$ and $U_2$ are
$G$-extensions of $U$, then $\Hom_H(U_1,U_2)\cong \R$.  Since $H$
acts trivially on $\Hom_H(U_1,U_2)$ and $G/H$ is of odd order,
$\Hom_H(U_1,U_2)$ must be a trivial $G$-module.  Therefore
$\Hom_G(U_1,U_2)$ is also isomorphic to $\R$, which means that
$U_1$ and $U_2$ are isomorphic as $G$-modules. This proves the
uniqueness of the $G$-extension of $U$.

It remains to prove~(2-b). Let $P$ be a normal subgroup of $G$
which contains $H$ and $P/H$ is a normal cyclic subgroup of $G/H$
of order $n$.  By~(1) above, $U$ has a $P$-extension, say $W$.
Then $\ind_P^GW$ is a $G$-extension of $2U$.
\end{proof}

\section{Fiber modules}

In this section we prove Theorem~\ref{theo:main_theorem} except
for Cases~A and~B.  The following two propositions can be proved
by the same argument as in the complex case, see Theorem~A and
its subsequent remark in~\cite{CKMS99a}.

\begin{prop} \label{prop:condition_of_fiber_$H$-module}
A real $H$-module is the fiber $H$-module of a real $G$-vector
bundle over $S(\rho)$ if and only if its character is
$G$-invariant and $G_1$-extendible (and $G_\mu$-extendible when
$\rho(G)=D_n$). \qed
\end{prop}

\begin{prop} \label{prop:construction}
Given $G_z$-extensions $V_z$ of a real $H$-module with
$G$-invariant character for $z=1$ (and $\mu$ when $\rho(G)=D_n$),
there exists a real $G$-vector bundle $E$ over $S(\rho)$ such that
the fiber $E_z$ of $E$ over $z$ is isomorphic to $V_z$ as
$G_z$-modules. \qed
\end{prop}

Proposition~\ref{prop:condition_of_fiber_$H$-module} gives a
characterization of the fiber $H$-module of a real $G$-vector
bundle over $S(\rho)$, and Proposition~\ref{prop:construction}
shows that the semi-group homomorphism $\Gamma$ in the
introduction is surjective. On the other hand, the following
proposition shows that $\Gamma$ is injective except for Cases~A
and~B, which proves Theorem~1.1 except for Cases~A and~B.

\begin{prop} \label{prop:real_isomorphism_theorem}
Let $\chi$ be a real irreducible character of $H$ which is
$G$-invariant. Except for Cases~A and~B, two real $G$-vector
bundles $E$ and $E'$ in $\Vect_G(S(\rho),\chi)$ are isomorphic if
and only if the fiber $G_z$-modules $E_z$ and $E'_z$ at $z\in
S(\rho)$ are isomorphic for $z=1$ (and for $z=\mu$ when
$\rho(G)=D_n$).
\end{prop}

\begin{proof}
The proof of~\cite[Theorem~6.1]{CKMS99a} holds in the real
category with slight modification. For reader's convenience we
shall give the argument when $\rho(G)$ is finite. The case when
$\rho(G)$ is infinite is easy since the action of $G$ on
$S(\rho)$ is transitive, see~\cite[Proposition~2.3]{CKMS99a} for
details.

We first note that if there exists an equivariant isomorphism
$\Psi\colon E\to E'$, then it must satisfy the equivariance
condition $\Psi_{\rho(g)z}=g\Psi_z g^{-1}$ for any $g\in G$ where
$\Psi_z=\Psi|_{E_z}$.

Suppose $\rho(G)\subsetneq SO(2)$. Then $G_1=H$, $\rho(G)$ is
finite cyclic, say, of order $n$, and since Case~A is excluded,
$\chi$ is not of real type. Choose an element $a\in G$ such that
$\rho(a)$ is the rotation through the angle $2\pi/n$.  By the
assumption we have an $H$-linear isomorphism $\Psi_1\colon E_1\to
E'_1$. Set $\Psi_{\rho(a)1}=a\Psi_1 a^{-1}$, which is also an
$H$-linear isomorphism. Then we connect $\Psi_1$ and
$\Psi_{\rho(a)1}$ continuously in the set of $H$-linear
isomorphisms of the fiber $H$-module along the arc of $S(\rho)$
joining $1$ and $\rho(a)1=e^{2\pi i/n}$. This is possible because
the set of $H$-linear isomorphisms of the fiber $H$-module is
homeomorphic to a general linear group over $\C$ or $\Q$
depending on the type of $\chi$ (remember that $\chi$ is not of
real type), and it is arcwise connected. Thus we have a bundle
isomorphism between $E$ and $E'$ restricted to the arc of
$S(\rho)$ joining $1$ and $\rho(a)1$. We extend this isomorphism
to an entire isomorphism over $S(\rho)$ using the equivariance
condition $\Psi_{\rho(a)z}=a\Psi_z a^{-1}$.

When $\rho(G)=D_n$, we choose a $G_1$-linear isomorphism $\Psi_1$
and a $G_\mu$-linear isomorphism $\Psi_\mu$. Similarly to the
above, we connect $\Psi_1$ and $\Psi_\mu$ as $H$-linear
isomorphisms along the arc of $S(\rho)$ joining $1$ and
$\mu=e^{\pi i/n}$, and then extend it to an isomorphism over
$S(\rho)$ using the equivariance condition. But it is not always
possible to connect $\Psi_1$ and $\Psi_\mu$ when $\chi$ is of
real type because the set of $H$-linear isomorphisms of the fiber
$H$-module, which is homeomorphic to $\GL(m,\R)$, is not arcwise
connected. In this case, however, we have another assumption that
$\chi$ is $G_z$-extendible for $z=1$ or $\mu$ since Case~B is
excluded. By Lemma~\ref{lemm:real_Mackey}~(1-a) $\chi$ has two
$G_z$-extensions, say $\widetilde\chi_1$ and $\widetilde\chi_2$.
Thus the character of $E_z\cong E'_z$ as a $G_z$-module is of the
form $m_1\widetilde\chi_1+m_2\widetilde\chi_2$ for some
nonnegative integers $m_1$ and $m_2$ with $m=m_1+m_2$, so that the
set of $G_z$-linear isomorphisms between $E_z$ and $E'_z$ is
homeomorphic to $\GL(m_1,\R)\times\GL(m_2,\R)$. Since the
inclusion map from $\GL(m_1,\R)\times\GL(m_2,\R)$ to $\GL(m,\R)$
induces a surjection on the $\pi_0$ level, it is possible to
choose a $G_z$-linear isomorphism $\Psi_z$ so that $\Psi_1$ and
$\Psi_\mu$ can be connected in the set of $H$-linear isomorphisms
of the fiber $H$-module.
\end{proof}

\section{Topological complexity: Cases~A and~B}

Propositions~\ref{prop:construction}
and~\ref{prop:real_isomorphism_theorem} show that the map
$\Gamma$ is an isomorphism except for Cases~A and~B. In this
section we investigate the structure on $\Vect_G(S(\rho),\chi)$
for Cases~A and~B, and complete the proof of
Theorem~\ref{theo:main_theorem}.

\subsection*{Case~A}
The case where $\rho(G)\subsetneq SO(2)$ and $\chi$ is of real type. In
this case one can reduce the study of $\Vect_G(S(\rho),\chi)$ to
the nonequivariant case.

\begin{lemm} \label{lemm:Case_A}
In Case~A, the semi-group $\Vect_G(S(\rho),\chi)$ is generated by
two elements $N_\chi^\pm$ with relation $2N_\chi^+=2N_\chi^-$.
Moreover, $N_\chi^\pm$ have $\chi$ as the character of the fiber
$H$-modules, and they are related in such a way that $N_\chi^-$
is obtained from $N_\chi^+$ by tensoring with a nontrivial real
$G$-line bundle over $S(\rho)$ with trivial fiber $H$-module.
\end{lemm}

\begin{proof}
There is an element $L$ in $\Vect_G(S(\rho),\chi)$ with $\chi$ as the
character of the fiber $H$-module by
Proposition~\ref{prop:condition_of_fiber_$H$-module},
and we have the semi-group isomorphisms
\[
\Vect_G(S(\rho),\chi) \cong \Vect_{G/H}(S(\rho)) \cong \Vect(S^1),
\]
where the former isomorphism is given by sending $E$ to
$\Hom_H(L,E)$ and the latter is given by taking orbit spaces by
the $G/H$-action. In fact, the map sending $F\in
\Vect_{G/H}(S(\rho))$ to $L\otimes F$ is the inverse of the former
isomorphism, where $F$ is viewed as a $G$-vector bundle through
the quotient map from $G$ to $G/H$ (see Lemma~2.2
in~\cite{CKMS99a} for details), and the latter is an isomorphism
because the action of $G/H$ on $S(\rho)$ is free. As is well
known, $\Vect(S^1)$ is generated by the trivial line bundle
$\epsilon$ and the Hopf line bundle $\eta$ with relation
$2\epsilon=2\eta$. Therefore, if we denote by $N_\chi^\pm$ the
two generators of $\Vect_G(S(\rho),\chi)$ corresponding to
$\epsilon$ and $\eta$ in $\Vect(S^1)$ through the above
isomorphism, then the lemma follows except the last statement. To
see the last statement, we note that $\Hom_H(N_\chi^+,N_\chi^-)$
is a nontrivial real $G$-line bundle over $S(\rho)$ with trivial
fiber $H$-module and that
\[
N_\chi^+ \otimes \Hom_H(N_\chi^+,N_\chi^-) \cong N_\chi^-,
\]
proving the last statement.
\end{proof}

\subsection*{Case~B}
The case where $\rho(G)=D_n$, $\chi$ is of real type, and neither
$G_1$- nor $G_\mu$-extendible. In this case we investigate
complex structures on the bundles in $\Vect_G(S(\rho),\chi)$.

Let $\F=\R$ or $\C$, and set
\[
\CS(\F^k)\equiv\{J\in\GL(k,\F) \mid J^2=-I\},
\]
which is the set of complex structures on $\F^k$. Needless to
say, $\CS(\R^k)$ is empty unless $k$ is even. Viewing $\C$ as
$\R^2$ in a natural way induces an injective homomorphism from
$\GL(k,\C)$ to $\GL(2k,\R)$, so that it induces an injection from
$\CS(\C^k)$ to $\CS(\R^{2k})$ and we view $\CS(\C^k)$ as a subset
of $\CS(\R^{2k})$ through this map.

\begin{lemm} \label{lemm:complex structure}
\begin{enumerate}
\item $\CS(\C^k)$ has $k+1$ connected components.
\item $\CS(\R^{2k})$ has two connected components.
\item If $k$ is odd, then each connected component of $\CS(\R^{2k})$
contains $(k+1)/2$ connected components of $\CS(\C^k)$, while
if $k$ is even, then one connected component of $\CS(\R^{2k})$ contains
$k/2$ and the other contains $k/2+1$ connected components of $\CS(\C^k)$.
\end{enumerate}
\end{lemm}

\begin{proof}
(1) We note that $\GL(k,\C)$ acts on $\CS(\C^k)$ by conjugation.
Since the minimal polynomial of any element in $\CS(\C^k)$ has
distinct root it is diagonalizable. So two elements in $\CS(\C^k)$
are in the same orbit if and only if they have the same
eigenvalues which are $\pm i$ because $J^2=-I$. This implies that
$\CS(\C^k)$ has exactly $k+1$ connected components because there
are $k+1$ possibilities of the $k$ eigenvalues.

(2) $\GL(2k,\R)$ acts transitively on $\CS(\R^{2k})$, and the
isotropy subgroup at an element of $\CS(\R^{2k})$ is isomorphic to
$\GL(k,\C)$; so $\CS(\R^{2k})$ is homeomorphic to a homogeneous
space $\GL(2k,\R)/\GL(k,\C)$ which has two connected components
(see~\cite[Proposition~2.48]{McSa99} for more details).

(3) As observed in~(1) above, $k+1$ elements
\[
\diag(i,i,\dotsc,i),\diag(-i,i,\dotsc,i),\dotsc,\diag(-i,-i,\dotsc,-i)
\]
respectively lie in the $k+1$ different connected components of
$\CS(\C^k)$. Through the inclusion map from $\CS(\C^k)$ to
$\CS(\R^{2k})$, they respectively are mapped to
\[
\diag(J_0,J_0,\dotsc,J_0),\diag(-J_0,J_0,\dotsc,J_0),\dotsc,
\diag(-J_0,-J_0,\dotsc,-J_0)
\]
where $J_0$ is the $2\times 2$ matrix
$\begin{pmatrix}0&-1\\1&0\end{pmatrix}$. Since $-J_0$ and $J_0$
are conjugate by $\begin{pmatrix}0&1\\1&0\end{pmatrix}$ whose
determinant is negative, the $k+1$ elements above in $\CS(\C^k)$
are in a same connected component of $\CS(\R^{2k})$ if and only
if the number of $J_0$'s as entries are congruent modulo $2$.
This implies~(3).
\end{proof}

For a real $G$-module $V$, we denote the set of $G$-invariant
complex structures on $V$ by
\[
\CS(V)^G\equiv\{J\in\GL(V)^G\mid J^2=-I\},
\]
where $\GL(V)^G$ denotes the $G$-linear automorphisms of $V$. A
pair $(V,J)$ is a complex $G$-module whose realification is $V$.
We note that $\GL(V)^G$ acts on $\CS(V)^G$ by conjugation and that
two complex $G$-modules $(V,J)$ and $(V,J')$ are isomorphic if and
only if $J$ and $J'$ are in the same orbit of the $\GL(V)^G$
action.

We consider the following setting for later use.

\begin{lemm} \label{lemm:equivariant complex structure}
Let $K$ be a group and let $H$ be a normal subgroup of $K$.
Suppose
\begin{enumerate} \renewcommand{\theenumi}{\alph{enumi}}
\item $W$ is an irreducible real $K$-module of complex type,
\item $U$ is an irreducible real $H$-module of real type,
\item $\res_HW\cong 2U$.
\end{enumerate}
Then, for any positive integer $k$, $\CS(kW)^K$ can naturally be viewed as
a subspace of $\CS(2kU)^H$, and we have
\begin{enumerate}
\item $\CS(kW)^K$ has $k+1$ connected components,
\item $\CS(2kU)^H$ has two connected components,
\item if $k$ is odd, then each connected component of $\CS(2kU)^H$
contains $(k+1)/2$ connected components of $\CS(kW)^K$, while if
$k$ is even, then one connected component of $\CS(2kU)^H$ contains
$k/2$ and the other contains $k/2+1$ connected components of
$\CS(kW)^K$.
\end{enumerate}
\end{lemm}

\begin{proof}
It follows from the assumptions (a) and (b) that $\GL(kW)^K\cong
\GL(k,\C)$ and $\GL(2kU)^H\cong \GL(2k,\R)$. Therefore the lemma
follows from Lemma~\ref{lemm:complex structure}.
\end{proof}

We return to the original setting of Case~B. Denote by $U$ a real
irreducible $H$-module with $\chi$ as its character. Since $\chi$
is neither $G_1$- nor $G_\mu$-extendible and of real type,
$\ind_H^{G_z}U$ is the unique $G_z$-extension of $2U$ of complex
type for $z=1$ and $\mu$ by Lemma~\ref{lemm:real_Mackey}~(2).
Therefore we are in a setting to which Lemma~\ref{lemm:equivariant
complex structure} can be applied. Moreover this shows that an
element $E$ in $\Vect_G(S(\rho),\chi)$ must have the fibers at
$z=1$ and $\mu$ isomorphic to $k(\ind_H^{G_z}U)$ for some integer
$k$. In particular, its fiber $H$-module is $2kU$.

A $G$-invariant complex structure on $E$ is a $G$-vector bundle
automorphism $J$ of $E$ such that $J^2=-I$. A pair $(E,J)$ is a
complex $G$-vector bundle whose realification is $E$. We say that
two $G$-invariant complex structures $J$ and $J'$ on $E$ are
equivalent if $(E,J)$ and $(E,J')$ are isomorphic as complex
$G$-vector bundles, in other words, if $J$ and $J'$ are conjugate
by a real $G$-vector bundle automorphism of $E$.

\begin{lemm} \label{lemm:number_of_complex_structures}
The number of inequivalent $G$-invariant complex structures on
$E$ is $(k+1)^2/2$ if $k$ is odd, and $k(k/2+1)$ or $k(k/2+1)+1$
if $k$ is even.
\end{lemm}

\begin{proof}
Let $\CS(E_z)$ be the set of (not necessarily invariant) complex
structures on the fiber $E_z$. The collection $\CS(E)$ of
$\CS(E_z)$ over $z\in S(\rho)$ forms a $G$-fiber bundle over
$S(\rho)$, the $G$-action on $\CS(E)$ being induced from that on
$E$. Then a $G$-invariant complex structure on $E$ can be viewed
as a continuous $G$-equivariant cross section of the $G$-fiber
bundle. The image of the cross section lies in $\CS(E)^H$ because
$H$ acts trivially on $S(\rho)$.

In order to construct a continuous $G$-equivariant cross section
of $\CS(E)\to S(\rho)$, we choose a pair of points from
$\CS(E_1)^{G_1}$ and $\CS(E_\mu)^{G_\mu}$ (i.e., one point from
each), which can be connected by a continuous cross section of
$\CS(E)^H$ restricted to the arc $R$ in $S(\rho)$ joining $1$ and
$\mu=e^{\pi i/n}$. Not all pairs of those points are connected by
such a cross section as observed later. But, once we find such a
cross section, we can extend it to an entire $G$-equivariant cross
section using the equivariance as is done in the proof of
Proposition~\ref{prop:real_isomorphism_theorem}. On the other
hand, we know in~\cite[Theorem~6.1]{CKMS99a} that isomorphism
classes of complex $G$-vector bundles over $S(\rho)$ are
distinguished by the complex fiber $G_1$- and $G_\mu$-modules.
Therefore, the number $CS(E)^G$ of inequivalent $G$-invariant
complex structures on $E$ is equal to the number of pairs of
connected components in $\CS(E_1)^{G_1}$ and $\CS(E_\mu)^{G_\mu}$
which are connected through $\CS(E)^H|_R$.

Suppose $k$ is even. Denote by $C_z^1$ and $C_z^2$, for $z=1$ and
$\mu$, the connected components of $\CS(E_z)^H|_R$ containing
$k/2$ and $k/2+1$ connected components of $\CS(E_z)^{G_z}$,
respectively. If $C_1^1$ and $C_\mu^1$ are connected through
$\CS(E)^H|_R$, then so are $C_1^2$ and $C_\mu^2$.  Counting the
number of choices of pairs of connected components in
$\CS(E_1)^{G_1}$ and $\CS(E_\mu)^{G_\mu}$ which are connected
through $\CS(E)^H|_R$, one has
\[
CS(E)^G=
\biggl(\frac{k}{2}\biggr)^2+\biggl(\frac{k}{2}+1\biggr)^2=k(k/2+1)+1.
\]
On the other hand, if $C_1^1$ and $C_\mu^2$ are connected through
$\CS(E)^H|_R$, then so are $C_1^2$ and $C_\mu^1$ and one has
\[
CS(E)^G=\frac{k}{2}\biggl(\frac{k}{2}+1\biggr)+\frac{k}{2}
\biggl(\frac{k}{2}+1\biggr)=k(k/2+1).
\]
For $k$ odd, a similar argument proves that
\[
CS(E)^G=\biggl(\frac{k+1}{2}\biggr)^2+\biggl(\frac{k+1}{2}\biggr)^2=(k+1)^2/2.
\]
\end{proof}

\begin{lemm} \label{lemm:Case_B}
In Case~B, the semi-group $\Vect_G(S(\rho),\chi)$ is generated by
two elements $M_\chi^\pm$ with relation $2M_\chi^+=2M_\chi^-$.
Moreover, $M_\chi^\pm$ have $2\chi$ as the character of the fiber
$H$-modules, and they are related in such a way that $M_\chi^-$
is obtained from $M_\chi^+$ by tensoring with a nontrivial
$G$-line bundle over $S(\rho)$ with trivial fiber $H$-module.
\end{lemm}

\begin{proof}
Since $\chi$ is of real type, the character of $U\otimes\C$ is
also $\chi$; so we may view $\chi$ as a complex irreducible
character of $H$. We have proved in~\cite[Theorem~B]{CKMS99a} that
the semi-group $\Vect_G^{\C}(S(\rho),\chi)$ of isomorphism
classes of complex $G$-vector bundles over $S(\rho)$ with
multiples of $\chi$ as the character of fiber $H$-modules is
generated by four elements $L_\chi^{\pm\pm}$ with relation
$L_\chi^{++}+L_\chi^{--}=L_\chi^{+-}+L_\chi^{-+}$, where
$L_\chi^{\pm\pm}$ are complex $G$-vector bundles over $S(\rho)$
with $U\otimes\C$ as the fiber $H$-module such that the fiber
$G_1$-modules (resp.~$G_\mu$-modules) of $L_\chi^{st}$ and
$L_\chi^{s't'}$, where $s$, $s'$, $t$, and $t'$ denote $+$ or $-$,
agree if and only if $s=s'$ (resp.~$t=t'$). In fact, the two
non-isomorphic fiber $G_1$-modules (resp.~$G_\mu$-modules) of
$L_\chi^{\pm\pm}$ are complex conjugate to each other, so the
complex conjugate (or~dual) bundles of $L_\chi^{++}$ and
$L_\chi^{+-}$ are respectively $L_\chi^{--}$ and $L_\chi^{-+}$.

Let $\Phi\colon\Vect_G^{\C}(S(\rho),\chi) \to
\Vect_G(S(\rho),\chi)$ be the realification map. It is surjective
by Lemma~\ref{lemm:number_of_complex_structures}. Since any
complex $G$-vector bundle is isomorphic to its complex conjugate
bundle as real $G$-vector bundles,
$\Phi(L_\chi^{++})=\Phi(L_\chi^{--})$ and
$\Phi(L_\chi^{+-})=\Phi(L_\chi^{-+})$. Therefore the relation
$L_\chi^{++}+L_\chi^{--}=L_\chi^{+-}+L_\chi^{-+}$ on
$\Vect_G^{\C}(S(\rho),\chi)$ reduces to
$2\Phi(L_\chi^{++})=2\Phi(L_\chi^{+-})$ on
$\Vect_G(S(\rho),\chi)$. It follows that for each fixed fiber
dimension there are at most two elements in
$\Vect_G(S(\rho),\chi)$. We claim that there is no other relation.
It suffices to show that there are exactly two elements in
$\Vect_G(S(\rho),\chi)$ for a fixed fiber dimension.  If there is
only one element for a fixed fiber dimension, say $2k\dim U$,
then the unique bundle must have $(k+1)^2$ inequivalent
$G$-invariant complex structures because the number of elements
in $\Vect_G^\C(S(\rho),\chi)$ of (real) fiber dimension $2k\dim
U$ is exactly $(k+1)^2$~\cite[Corollary~5.2]{CKMS99a}. This
contradicts Lemma~\ref{lemm:number_of_complex_structures}.

It remains to show that the two generators $\Phi(L_\chi^{++})$ and
$\Phi(L_\chi^{+-})$ are related by tensoring with a nontrivial
$G$-line bundle with trivial fiber $H$-module. The fiber
$G_1$-modules of $L_\chi^{+-}$ and $L_\chi^{++}$ at $1$ are
isomorphic but the fiber $G_\mu$-modules of them at $\mu$ are not,
more precisely, they are related through the tensor product with
the nontrivial real $1$-dimensional $G_\mu$-module defined by
$G_\mu\to G_\mu/H\cong\{\pm 1\}$, see
Lemma~\ref{lemm:real_Mackey}~(1). Therefore $\Phi(L_\chi^{+-})$
is obtained from $\Phi(L_\chi^{++})$ by tensoring with a real
$G$-line bundle with trivial fiber $H$-module, whose fiber at $1$
is the trivial $G_1$-module and the fiber at $\mu$ is the
nontrivial $G_\mu$-module. The existence of such line bundle is
guaranteed by Proposition~\ref{prop:construction}.
\end{proof}

\begin{proof}[Proof of Theorem~\ref{theo:main_theorem}]
The map $\Gamma$ is surjective by
Proposition~\ref{prop:construction} and injective except for
Cases~A and~B by Proposition~\ref{prop:real_isomorphism_theorem}.
In both Cases~A and~B the target of $\Gamma$ is a semi-group
generated by one element by Lemmas~\ref{lemm:structure_Rep(K)} and
\ref{lemm:structure_Rep(K1,K2)} while the domain of $\Gamma$ is
generated by two elements with the relation as in
Lemmas~\ref{lemm:Case_A} and~\ref{lemm:Case_B}. This implies that
$\Gamma$ is two to one.

Finally, we note that tensoring elements in
$\Vect_G(S(\rho),\chi)$ with a nontrivial $G$-line bundle with
trivial $H$-module does not change the fiber $G_1$-modules
(resp.~fiber $G_1$- and $G_\mu$-modules) in Case~A (resp.~Case~B).
This implies the last statement in the theorem.
\end{proof}

\section{The semi-group structure on $\Vect_G(S(\rho),\chi)$}

In this section we determine the semi-group structure on
$\Vect_G(S(\rho),\chi)$. Let $e_1$ and $e_\mu$ denote the numbers
of $G_1$- and $G_\mu$-extensions of $\chi$, respectively.  When
$\rho(G)$ agrees with $O(2)$ or is contained in $SO(2)$, we
define $e_\mu$ to be $1$ for convenience. In both real and
complex category, the semi-group structure on the target of
$\Gamma$ is determined by the numbers $e_1$ and $e_\mu$. The
numbers $e_1$ and $e_\mu$ depend only on the types of $\rho(G)$
in the complex category, but this is not true in the real
category. This is another complexity in our study arising from
real representation theory.

The possible values of $e_1$ and $e_\mu$ according to $\rho(G)$
and the type of $\chi$ are given by
Table~\ref{table:number_of_extensions}
and~\ref{table:number_of_extensions_by_type}.

\begin{table}[hbt]
\begin{tabular}{|c||c|c|c|c|c|c|c|c|c|} \hline
$(e_1,e_\mu)$ & $(0,0)$ & $(1,0)$ & $(0,1)$ & $(2,0)$ & $(0,2)$ &
$(1,1)$ & $(2,1)$ & $(1,2)$ & $(2,2)$ \\\hline\hline
$\rho(G)\subset SO(2)$ & $\times$ & $\times$ & $\times$ & $\times$ & $\times$ &
{\scriptsize$\bigcirc$} & $\times$ & $\times$ & $\times$ \\ \hline
$\rho(G)=O(2)$ & $\times$ & $\times$ & {\scriptsize$\bigcirc$} &
$\times$ & $\times$ & {\scriptsize$\bigcirc$} &
{\scriptsize$\bigcirc$} & $\times$ & $\times$ \\ \hline
$\rho(G)=D_n$ & {\scriptsize$\bigcirc$} & {\scriptsize$\bigcirc$}
& {\scriptsize$\bigcirc$} & {\scriptsize$\bigcirc$} &
{\scriptsize$\bigcirc$} & {\scriptsize$\bigcirc$} &
{\scriptsize$\bigcirc$} & {\scriptsize$\bigcirc$} &
{\scriptsize$\bigcirc$} \\\hline
\end{tabular} \vspace{1.25ex}
\caption{The possible values of $(e_1,e_\mu)$ according to $\rho(G)$}
\label{table:number_of_extensions}
\end{table}

\begin{table}[htb]
\begin{tabular}{|c||c|c|c|c|c|c|c|c|c|} \hline
$(e_1,e_\mu)$ & $(0,0)$ & $(1,0)$ & $(0,1)$ & $(2,0)$ & $(0,2)$ &
$(1,1)$ & $(2,1)$ & $(1,2)$ & $(2,2)$ \\\hline\hline real &
{\scriptsize$\bigcirc$} & $\times$ & $\times$ &
{\scriptsize$\bigcirc$} & {\scriptsize$\bigcirc$} & $\times$ &
$\times$ & $\times$ & {\scriptsize$\bigcirc$}
\\\hline
complex & {\scriptsize$\bigcirc$} & {\scriptsize$\bigcirc$} &
{\scriptsize$\bigcirc$} & {\scriptsize$\bigcirc$} &
{\scriptsize$\bigcirc$} & {\scriptsize$\bigcirc$} &
{\scriptsize$\bigcirc$} & {\scriptsize$\bigcirc$} &
{\scriptsize$\bigcirc$} \\\hline quaternionic & $\times$ &
$\times$ & $\times$ & $\times$ & $\times$ &
{\scriptsize$\bigcirc$} & {\scriptsize$\bigcirc$} &
{\scriptsize$\bigcirc$} & {\scriptsize$\bigcirc$} \\\hline
\end{tabular} \vspace{1.25ex}
\caption{The possible values of $(e_1,e_\mu)$ according to the type of
$\chi$ when $\rho(G)=D_n$}
\label{table:number_of_extensions_by_type}
\end{table}

We state here the semi-group structure on $\Vect_G(S(\rho),\chi)$
according to the values of $e_1$ and $e_\mu$.

\begin{theo} \label{theo:structure}
Except for Cases~A and~B, the semi-group $\Vect_G(S(\rho),\chi)$
is generated by
\begin{enumerate}
\item one element $L_\chi$, if $(e_1,e_\mu)=(0,0),(1,0),(0,1)$ or $(1,1)$,
\item two elements $L_\chi^\pm$, if $(e_1,e_\mu)=(2,1)$ or $(1,2)$,
\item three elements $\widetilde L_\chi^0,\widetilde L_\chi^\pm$
with relation $2\widetilde L_\chi^0=\widetilde L_\chi^++\widetilde
L_\chi^-$, if $(e_1,e_\mu)=(2,0)$ or $(0,2)$,
\item four elements $L_\chi^{\pm\pm}$ with relation
$L_\chi^{++}+L_\chi^{--}=L_\chi^{+-}+L_\chi^{-+}$, if
$(e_1,e_\mu)=(2,2)$.
\end{enumerate}
In Case~A or~B, it is generated by
\begin{enumerate} \setcounter{enumi}{4}
\item two elements $\widetilde{\widetilde L}{}_\chi^\pm$ with relation
$2\widetilde{\widetilde L}{}_\chi^+=2\widetilde{\widetilde
L}{}_\chi^-$.
\end{enumerate}
Moreover, except for $\widetilde L_\chi^0$ in the case~(3), all
generators are related through tensor product with real $G$-line
bundles over $S(\rho)$ with trivial fiber $H$-module.
\end{theo}

\begin{proof}
The statements (1)--(4) follow from
Lemmas~\ref{lemm:structure_Rep(K)},
\ref{lemm:structure_Rep(K1,K2)} and
Theorem~\ref{theo:main_theorem}, and the statement~(5) follows
from Lemmas~\ref{lemm:Case_A} and~\ref{lemm:Case_B}.

We now prove the last statement in the theorem. After setting
$K_1=G_1$ and $K_2=G_\mu$, it is obvious that the inverse images
of the pairs $(\R_1^\pm,\R_2^\pm)$ in the remark after
Lemma~\ref{lemm:structure_Rep(K1,K2)} by the semi-group
homomorphism $\Gamma$ in Theorem~\ref{theo:main_theorem} are real
$G$-line bundles with trivial fiber $H$-module. Moreover,
$\Gamma$ preserves the two tensor product operations, one on
$\Vect_G(S(\rho),\chi)$ with real $G$-line bundles and the other
on $\Rep(G_1,G_\mu,\chi)$ by the pairs $(\R_1^\pm,\R_2^\pm)$.
Therefore, Proposition~\ref{prop:real_isomorphism_theorem}
implies that, except for Cases~A and~B, the generators of
$\Vect_G(S(\rho),\chi)$ are related through tensor product with
real $G$-line bundles with trivial fiber $H$-module. The same
argument also holds for $\Rep(G_1,\chi)$ by
Lemma~\ref{lemm:structure_Rep(K)}. For Cases~A and~B, the
statement follows from the last statement of
Lemmas~\ref{lemm:Case_A} and~\ref{lemm:Case_B}.
\end{proof}

\begin{coro} \label{coro:enumeration}
Let $N$ be the number of isomorphism classes of real $G$-vector
bundles over $S(\rho)$ with $m\chi$ as the character of the fiber
$H$-modules. In case that $m$ is odd, $N$ is zero if $e_1=0$ or
$e_\mu=0$. Otherwise, except for Cases~A and~B,
\[
N = \begin{cases}
1, & \text{if $(e_1,e_\mu)=(0,0),(1,0),(0,1)$, or $(1,1)$}, \\
m+1, & \text{if $(e_1,e_\mu)=(2,0),(0,2),(2,1)$, or $(1,2)$}, \\
(m+1)^2, & \text{if $(e_1,e_\mu)=(2,2)$}.
\end{cases}
\]
In Case~A or~B, the number $N$ is exactly two.
\end{coro}

\begin{proof} The proof is elementary and left to the reader.
\end{proof}

\section{Triviality of real $G$-vector bundles over a circle}

In this section we investigate triviality of the generators in
Theorem~\ref{theo:structure} when $\rho(G)$ is finite. Triviality
of a $G$-vector bundle is closely related to the existence of a
$G$-extension of the fiber $H$-module in the following sense: For
a given $H$-module $V$, there exists at least one trivial
$G$-vector bundle with $V$ as its fiber $H$-module if $V$ extends
to a $G$-module. In the following we denote by $Z_n$ the finite
cyclic subgroup of $SO(2)$ generated by the rotation through an
angle $2\pi/n$. Then $\rho(G)=Z_n$ for some $n$ if
$\rho(G)\subsetneq SO(2)$. Denote by $1$ the trivial real
$H$-module of dimension one, in other words, $H$ acts trivially
on $1$. In the notation of Lemma~\ref{lemm:Case_A} and
Theorem~\ref{theo:structure}, real $G$-line bundles over
$S(\rho)$ with trivial fiber $H$-module are denoted by $N_1^\pm$
and $L_1^{\pm\pm}$ according as $\rho(G)=Z_n$ and $D_n$,
respectively.

\begin{lemm} \label{lemm:line_bundles}
\textup{(1)} Suppose $\rho(G)=Z_n$. If $n$ is even, then $N_1^\pm$
are both trivial. If $n$ is odd, then one of them, say $N_1^+$,
is trivial and $N_1^-$ is nontrivial.

\textup{(2)} Suppose $\rho(G)=D_n$.  If $n$ is even, then
$L_1^{\pm\pm}$ are all trivial. If $n$ is odd, then two of them,
say $L_1^{++}$ and $L_1^{--}$, are trivial and the other two are
nontrivial.
\end{lemm}

\begin{proof}
(1) Since $G/H$ acts freely on $S(\rho)$, every real $G$-line
bundle over $S(\rho)$ with trivial fiber $H$-module is the
pull-back of a real line bundle over $S^1$ by the quotient map
$\pi\colon S(\rho)\to S(\rho)/G\cong S^1$. Suppose $n$ is even.
Then $\pi^*\colon H^1(S(\rho)/G,\Z/2)\to H^1(S(\rho),\Z/2)$ is
trivial, so pullback line bundles by $\pi$ have trivial first
Whitney classes, which means that the underlying line bundles
over $S(\rho)$ are trivial. According
to~\cite[Proposition~1.1]{KiMa94}, an equivariant line bundle is
trivial if and only if its underlying bundle is trivial. Thus,
$N_1^\pm$ are both trivial when $n$ is even.

If $n$ is odd, then $\pi^*$ above is an isomorphism.  Therefore,
exactly one of $N_1^\pm$ has trivial first Whitney class. This
together with the result in~\cite{KiMa94} mentioned above shows
that exactly one of $N_1^\pm$ is trivial equivariantly.

(2) Set $P=\rho^{-1}(Z_n)$. Since $P/H$ acts freely on $S(\rho)$,
$L_1^{\pm\pm}$ are pullback of real $G/P$-line bundles over
$S(\rho)/P$ by the quotient map $\pi\colon S(\rho)\to
S(\rho)/P$.  Here $G/P$ is of order two and acts on the circle
$S(\rho)/P$ as reflection, so we may think of $G/P$ as $D_1$.
According to Theorem~\ref{theo:structure} (or
Corollary~\ref{coro:enumeration}) there are four real $D_1$-line
bundles over $S(\rho)/P$.  Since the map $\Gamma$ is an
isomorphism in this case, they are distinguished by their fiber
$D_1$-modules over the points $\pm 1\in S(\rho)/P$. More
precisely, there are two possibilities for the fiber
$D_1$-modules at $1$ and $-1$ respectively since there are two
real one-dimensional $D_1$-modules (the trivial one and the
nontrivial one), and hence altogether there are four real
$D_1$-line bundles over $S(\rho)/P$. Moreover, $D_1$-line bundles
are trivial if and only if the fiber $D_1$-modules at $\pm 1$ are
isomorphic (see also \cite{Kim94}).

If $n$ is even, then all pullback line bundles by $\pi$ are
trivial as discussed in~(1); so $L_1^{\pm\pm}$ are all trivial.
If $n$ is odd, then the pullback by $\pi$ preserves the
triviality of line bundles because $\pi^*\colon
H^1(S(\rho)/P;\Z/2)\to H^1(S(\rho);\Z/2)$ is an isomorphism.
Since there are exactly two trivial $D_1$-line bundles over
$S(\rho)/P$, two of $L_1^{\pm\pm}$ are trivial and the other two
are nontrivial.
\end{proof}

\begin{rema}
Suppose $\rho(G)=D_n$. For $z=1$ and $\mu$, denote by $\R_z^+$ and
$\R_z^-$, respectively, the trivial and the nontrivial real
$G_z$-module of dimension one with trivial $H$-action, see also
the remark after Lemma~\ref{lemm:structure_Rep(K1,K2)}. Then we
may assume without loss of generality that the images of
$L_1^{st}$ by $\Gamma$ in Theorem~\ref{theo:main_theorem} are
$(\R_1^s,\R_\mu^t)$, where $s$ and $t$ denote a sign $+$ or $-$.
\end{rema}

\begin{theo} \label{theo:triviality_for_L}
Let $\rho(G)=Z_n$ or $D_n$, and let $\chi$ be a real irreducible
character of $H$ which is $G$-invariant. If $n$ is even, then the
generators in Theorem~\ref{theo:structure} except for $\widetilde
L_\chi^0$ are all trivial or all nontrivial in each case. If $n$
is odd, then
\begin{enumerate}
\item $L_\chi$ is trivial,
\item $L_\chi^{\pm}$ are both trivial,
\item $\widetilde L_\chi^0$ and $\widetilde L_\chi^{\pm}$ are all trivial,
\item two of $L_\chi^{\pm\pm}$ are trivial and the other two are nontrivial,
\item one of $\widetilde{\widetilde L}{}_\chi^\pm$ is trivial and the
other is nontrivial.
\end{enumerate}
\end{theo}

\begin{proof}
Recall from the last statement in Theorem~\ref{theo:structure}
that all generators are related through tensor product with the
real $G$-line bundles $N_1^\pm$ and $L_1^{\pm\pm}$ according as
$\rho(G)=Z_n$ and $D_n$, respectively. These line
bundles are all trivial if $n$ is even by
Lemma~\ref{lemm:line_bundles}. So the existence of one trivial
generator implies triviality of the other generators, and this
finishes the proof in case that $n$ is even.

In the following we assume that $n$ is odd. Denote by $U$ a real
irreducible $H$-module with $\chi$ as its character. Recall that
the fiber $H$-module of a generator is $U$ if both $e_1$ and
$e_\mu$ are nonzero, and $2U$ otherwise. In case $\rho(G)=D_n$,
we choose elements $a$ and $b$ in $G$ such that $\rho(a)$ is the
rotation through the angle $2\pi/n$ and $\rho(b)$ is the
reflection about the $x$-axis. Then $G_1$ (resp.~$G_\mu$) is
generated by $H$ and $b$ (resp.~$ab$).

(1) It suffices to show that the fiber $H$-module of a generator
extends to a $G$-module. In case that $e_1=e_\mu=1$, the fiber
$H$-module of $L_\chi$ is $U$ and it is $G$-extendible by
Lemma~\ref{lemm:extension}. The other case is that either $e_1=0$
or $e_\mu=0$ and in this case the fiber $H$-module of $L_\chi$ is
$2U$ which is $G$-extendible by Lemma~\ref{lemm:extension}~(2-b).

(2) In this case $\rho(G)=D_n$ by
Table~\ref{table:number_of_extensions} and the fiber $H$-modules
of generators are $U$ which is $G$-extendible by
Lemma~\ref{lemm:extension}~(2). So there is at least one trivial
generator, say $L_\chi^+$. Since $(e_1,e_\mu)=(2,1)$ or $(1,2)$,
the tensor product of $L_\chi^+$ with $L_1^{--}$ has different
fiber $G_z$-module from that of $L_\chi^+$ at the point $z$ such
that $e_z=2$. So we get the other generator $L_\chi^-\cong
L_\chi^+\otimes L_1^{--}$. Since $L_1^{--}$ is trivial by
Lemma~\ref{lemm:line_bundles}, so is $L_\chi^-$.

(3) In this case $\rho(G)=D_n$ by
Table~\ref{table:number_of_extensions} and the fiber $H$-modules
of generators are $2U$ because either $e_1$ or $e_\mu$ is zero.
Set $P=\rho^{-1}(Z_n)$. Then $U$ has a $P$-extension, say $V$, by
Lemma~\ref{lemm:extension}~(1). Note that the fiber modules of
$\widetilde L_\chi^0$ at $1$ and $\mu$ are isomorphic to
$\ind_H^{G_1}U$ and $\ind_H^{G_\mu}U$, respectively. Thus
$\widetilde L_\chi^0$ is isomorphic to the product bundle
$S(\rho)\times \ind_P^G V$ by
Proposition~\ref{prop:real_isomorphism_theorem}.

We next consider triviality of the generators $\widetilde
L_\chi^\pm$. It suffices to show that at least one generator, say
$\widetilde L_\chi^+$, is trivial. Then so is the other generator
$\widetilde L_\chi^-\cong \widetilde L_\chi^+\otimes L_1^{--}$.
We assume that $(e_1,e_\mu)=(2,0)$. The other case $(e_1,e_\mu)=
(0,2)$ can be proved similarly.

\begin{clai}
$\chi$ is of real type.
\end{clai}

\begin{proof}[Proof of Claim]
Since $\chi$ is not of quaternionic type by
Table~\ref{table:number_of_extensions_by_type}, it suffices to
prove that $\chi$ is not of complex type. Suppose that $\chi$ is
of complex type.  Then there is a complex $H$-module $V$ such
that $U\otimes\C\cong V\oplus\overline V$ and $V\ncong\overline
V$ as complex $H$-modules. We note that the realifications of $V$
and $\overline V$ are $U$, and since $\chi$ is $G$-invariant,
$\side{g}\chi_V=\chi_V$ or $\chi_{\overline V}$ for $g\in G$
where $\chi_V$ and $\chi_{\overline V}$ denote the characters of
$V$ and $\overline V$ respectively.

Since $e_1=2$ and $e_\mu=0$ by assumption, $U$ has two
$G_1$-extensions of complex type by Lemma~\ref{lemm:real_Mackey}
but no $G_\mu$-extension. It follows that $V$ is $G_1$-extendible
but not $G_\mu$-extendible, so $\chi_V$ is $G_1$-invariant but
not $G_\mu$-invariant. Namely, $\side{b}\chi_V=\chi_V$ and
$\side{ab}\chi_V=\chi_{\overline V}$, so that
$\side{a}\chi_V=\chi_{\overline V}$. Therefore
$\siden{a^n}\chi_V=\chi_{\overline V}$ because $n$ is odd. On the
other hand, since $a^n$ is an element of $H$,
$\siden{a^n}\chi_V=\chi_V$.  Therefore $\chi_V=\chi_{\overline
V}$, but this contradicts that $V\ncong \overline V$.  Thus
$\chi$ must be of real type.
\end{proof}

Since $\chi$ is of real type by the claim above, $U\otimes\C$ is
irreducible and its character is $G$-invariant. It follows that
there is a trivial complex $G$-vector bundle $F$ over $S(\rho)$
with $U\otimes\C$ as the fiber $H$-module,
see~\cite[Theorem~C~(3)]{CKMS99a}. Since $e_1=2$, there are two
$G_1$-extensions of $U$, say $\widetilde U_1$ and $\widetilde
U_2$. Their complexifications $\widetilde U_1\otimes\C$ and
$\widetilde U_2\otimes\C$ are non-isomorphic because $\widetilde
U_1\ncong\widetilde U_2$. Moreover these modules are both
$G_1$-extensions of $U\otimes\C$. Thus the fiber $G_1$-module,
say $F_1$, of $F$ at $1$ must be either $\widetilde U_1\otimes\C$
or $\widetilde U_2\otimes\C$. It follows that the realification
of $F_1$ is either $2\widetilde U_1$ or $2\widetilde U_2$.
Therefore the realification of $F$, which is trivial, is
isomorphic to one of $\widetilde L_\chi^\pm$.

(4) In this case $\rho(G)=D_n$ by
Table~\ref{table:number_of_extensions} and the fiber $H$-modules
of the generators are $U$. By a similar argument to the case~(2)
there are two trivial generators, $L_\chi^{++}$ and
$L_\chi^{--}\cong L_\chi^{++}\otimes L_1^{--}$. It suffices to
show that the other two generators $L_\chi^{+-}\cong
L_\chi^{++}\otimes L_1^{+-}$ and $L_\chi^{-+}\cong
L_\chi^{++}\otimes L_1^{-+}$ are nontrivial. Consider the
following isomorphisms
\begin{equation}
\Hom_H(L_\chi^{++},L_\chi^{+-})
\cong\Hom_H(L_\chi^{++},L_\chi^{++}\otimes L_1^{+-})
\cong\Hom_H(L_\chi^{++},L_\chi^{++})\otimes L_1^{+-}. \tag{**}
\end{equation}
Note that $\Hom_H(L_\chi^{++},L_\chi^{++})$ is isomorphic to the
product bundle $S(\rho)\times\R^k$, where $k=1,2$, or $4$
according to the type of $\chi$. It follows that
$\Hom_H(L_\chi^{++},L_\chi^{+-})\cong kL_1^{+-}$.

\begin{clai} \label{clai:nontriviality_of_L_1^+-}
Both $kL_1^{+-}$ and $kL_1^{-+}$ are nontrivial for all $k>0$.
\end{clai}

\begin{proof}[Proof of Claim]
Note that the fiber $G_1$-module of $L_1^{+-}$ at $1\in S(\rho)$
is the trivial $G_1$-module $\R_1^+$ while the fiber
$G_\mu$-module at $\mu$ is the nontrivial $G_\mu$-module
$\R_\mu^-$ by the remark after Lemma~\ref{lemm:line_bundles}.
Then $b$ (resp.~$ab$) acts on $\R_1^+$ (resp.~$\R_\mu^-$) as
multiplication by $1$ (resp.~$-1$). Recall that $H$ acts on both
$\R_1^+$ and $\R_\mu^-$ trivially, i.e., as multiplication by $1$.

Assume that $kL_1^{+-}$ is trivial. Then there exists a
$G$-module $W$ such that $\res_{G_1}W\cong k\R_1^+$ and
$\res_{G_\mu}W\cong k\R_\mu^-$. Thus $b$ and $ab$ act on $W$ as
multiplication by $1$ and $-1$, respectively. Hence $a$ acts on
$W$ as multiplication by $-1$, and since $n$ is odd, $a^n$ also
acts on $W$ as multiplication by $-1$. But this contradicts that
$a^n\in H$ acts trivially on $W$. In the same way we can prove
that $kL_1^{-+}$ is also nontrivial.
\end{proof}

Since $L_\chi^{++}$ is trivial, $L_\chi^{+-}$ must be nontrivial
by the equation~(**) and the claim above. Replacing $L_\chi^{+-}$
by $L_\chi^{-+}$ we can similarly prove  that $L_\chi^{-+}$ is
nontrivial.

(5) \emph{Case A:} In this case $\widetilde{\widetilde
L}{}_\chi^\pm$ are $N_\chi^\pm$ in Lemma~\ref{lemm:Case_A}. Since
$n$ is odd, $U$ has a $G$-extension by
Lemma~\ref{lemm:extension}~(1). So we may assume that one
generator, say $N_\chi^+$, is trivial. Then the following
isomorphisms
\[
\Hom_H(N_\chi^+,N_\chi^-) \cong \Hom_H(N_\chi^+,N_\chi^+\otimes
N_1^-) \cong \Hom_H(N_\chi^+,N_\chi^+)\otimes N_1^- \cong N_1^-
\]
imply that $N_\chi^-$ is nontrivial since $N_1^-$ is nontrivial by
Lemma~\ref{lemm:line_bundles}~(1).

\emph{Case B:} In this case $\widetilde{\widetilde L}{}_\chi^\pm$
are $M_\chi^\pm$ in Lemma~\ref{lemm:Case_B}. Remember that
$M_\chi^+=\Phi(L_\chi^{++})=\Phi(L_\chi^{--})$ and
$M_\chi^-=\Phi(L_\chi^{+-})=\Phi(L_\chi^{-+})$, see the proof of
Lemma~\ref{lemm:Case_B}. Since
$L_\chi^{++}\in\Vect_G^\C(S(\rho),\chi)$ is trivial by Theorem~C
in~\cite{CKMS99a}, $M_\chi^+$ is also trivial. In the following
we shall prove that $M_\chi^-$ is nontrivial.

Assume that $M_\chi^-$ is trivial, i.e., it is isomorphic to the
product bundle $S(\rho)\times W$ for some $G$-extension $W$ of the
fiber $H$-module $2U$.

\begin{clai} $W$ is of real type.
\end{clai}

\begin{proof}[Proof of Claim]
If $W$ is not of real type, then we may view $M_\chi^-$ as the
realification of a complex product bundle $S(\rho)\times W$, but
this contradicts that $M_\chi^-$ is the realification of the
nontrivial bundles $L_\chi^{+-}$ and $L_\chi^{-+}$ in
$\Vect_G^\C(S(\rho),\chi)$.
\end{proof}

Denote by $\chi_W$ the character of $W$. Every fiber $G_z$-module
of $M_\chi^-$, which is $\res_{G_z}W$, is isomorphic to
$\ind_H^{G_z}U$ and it is irreducible of complex type by
Lemma~\ref{lemm:real_Mackey}~(2).  It is well known in
representation theory that the character of
$\res_{G_z}W\cong\ind_H^{G_z}U$ is zero on $G_z\setminus H$. Thus
$\chi_W$ is always zero on $\bigcup_{z\in S(\rho)}G_z\setminus
H=G\setminus P$, where $P=\rho^{-1}(Z_n)$. It follows that we have
\[
1=\int_G \chi_W(g)^2 dg = \frac{1}{2} \int_P \chi_W(p)^2 dp +
\frac{1}{2} \int_{G\setminus P} \chi_W(p)^2 dp = \frac{1}{2}
\int_P \chi_W(p)^2 dp,
\]
so $\displaystyle{\int_P\chi_W(p)^2 dp=2}$. This implies
that $\res_PW$ is either irreducible of complex type or
reducible with different direct summands of real type.
In the sequel we show that neither case occurs.

It is easy to see that the latter case does not occur because if
it does, then each summand of $\res_PW$ is a $P$-extension of $U$
which contradicts the uniqueness of the $P$-extension of $U$ by
Lemma~\ref{lemm:extension}~(1).

Now, suppose $\res_PW$ is irreducible and of complex type. We
claim that the set $\CS(W)^G$ of $G$-invariant complex structures
on $W$ is not empty. Then it contradicts that $W$ is of real type.
Since $\res_{G_1}W$ is irreducible and of complex type, there
exists a $G_1$-invariant complex structure on $W$, i.e.,
$(\CS(W)^H)^{G_1/H}=\CS(W)^{G_1}\neq\varnothing$. This means that
each connected component of $\CS(W)^H\cong\CS(\R^2)$ is invariant
under the action of $G_1$ because $\CS(W)^H\cong\CS(\R^2)$ has
two connected components. On the other hand, since the order $n$
of $P/H$ is odd and the number of connected components of
$\CS(W)^H$ is two, each connected component is also invariant
under the $P/H$-action. Therefore, it is invariant under the
$G/H$-action because $P$ and $G_1$ generate $G$. Now we note that
each connected component of $\CS(W)^H\cong\CS(\R^2)$ is
homeomorphic to $\R^2$ (see~\cite[Exercise~2.57]{McSa99}) and that
any smooth action of a finite group on $\R^2$ is linear, so the
$G/H$-action on $\CS(W)^H$ has a fixed point, i.e.,
$\CS(W)^G\neq\varnothing$.
\end{proof}

\providecommand{\bysame}{\leavevmode\hbox to3em{\hrulefill}\thinspace}

\end{document}